%------------------------------------------------------------------------------
% ladder.tex
%                    
% Beginn: 13.02.2007
%------------------------------------------------------------------------------
% mredondo@criba.edu.ar  Maria Julia
\magnification=\magstep1   
\input amstex
\UseAMSsymbols
\input pictex
%\input german 
%\hoffset=0truecm \voffset=0truecm 
%\vsize=23truecm
\NoBlackBoxes
\parindent=18pt
\font\kl=cmr8

\font\gross=cmcsc10 scaled\magstep2
\font\abs=cmcsc10

\def\Hom{\operatorname{Hom}}
\def\End{\operatorname{End}}
\def\Ext{\operatorname{Ext}}
\def\Im{\operatorname{Im}}

\def\arr#1#2{\arrow <1.5mm> [0.25,0.75] from #1 to #2}

\def\t#1{\quad\text{#1}\quad} 
    %Fuer Text innerhalb von $$...$$. Verwendung: \T{oder}
   
          %Zum Einruecken, Verwendung: \E{dieses wird eingerueckt}
\vglue4truecm
\centerline{\gross The Ladder Construction of Pr\"ufer Modules.}     
\bigskip\medskip
\centerline{\abs Claus Michael Ringel}    
 \bigskip
\centerline{Dedicated to Maria Inez Platzeck on the occasion of her 60th birthday}
	\bigskip\bigskip
{\narrower\narrower\narrower \kl Abstract. Let $\ssize R$ be a ring
(associative, with 1). A non-zero 
module $\ssize M$ is said to
be a Pr\"ufer module provided there exists a surjective, locally nilpotent
endomorphism with kernel of finite length. The aim of this note is construct
Pr\"ufer modules starting from a pair of module homomorphisms 
$\ssize w,v\: U_0 \to U_1$,
where $\ssize w$ is injective and its cokernel is of finite length.
For $\ssize R = \Bbb Z$ the ring of integers, one can construct in this way
the ordinary Pr\"ufer groups considered in abelian group theory. Our
interest lies in the case that $\ssize R$ is an artin algebra. 
\par}
	\bigskip\bigskip
%------------------------------------------------------------------------------
	\bigskip
{\bf 1. The construction.}
	\medskip
Let $R$ be a ring (associative, with 1). The modules to be considered will usually
be left $R$-modules.
Our main interest will be the case where
$R$ is an artin algebra, however the basic construction should be of interest for any
ring $R$. In fact, the standard examples of what we call Pr\"ufer modules are the 
Pr\"ufer groups in abelian group theory, thus $\Bbb Z$-modules. Here is the definition
of a Pr\"ufer module: it is a non-zero 
module $P$ which has a surjective, locally nilpotent
endomorphism $\phi$ with kernel of finite length. If $H$ is the kernel of $\phi$, we 
often will write $P = H[\infty]$, and we will denote the kernel of $\phi^t$ by $H[t]$. 
Observe the slight ambiguity: given a Pr\"ufer module $P$, not only $\phi$ but also all
non-trivial powers of $\phi$ and maybe many other endomorphisms will have the required 
properties (surjectivity, locally nilpotency, finite length kernel).
	
The content of the paper is as follows. In the first section we show that any pair of
module homomorphisms $w,v:U_0 \to U_1$, where $w$ is injective with non-zero cokernel
of finite length, gives rise to a Pr\"ufer module. Section 2 provides some examples and
section 3 outlines the relationship between Pr\"ufer modules and 
various sorts of self-extensions of finite length modules. The final sections 4 and 5 
deal with degenerations in the sense of Riedtmann-Zwara: we will show that this degeneration
theory is intimately connected to the existence of Pr\"ufer modules with some splitting
property, and we will exhibit an extension of a recent result by Bautista
and Perez. Our interest in the questions considered here was stimulated by a series of lectures 
by Sverre Smal{\o} [S] at the Mar del Plata conference, March 2006, and we are indebted to
him as well as to M.C.R.Butler and G.Zwara for helpful comments.

For the relevance
of Pr\"ufer modules when dealing with artin algebras of infinite representation
type, we refer to a forthcoming paper [R5]. The appendix to section 3.3
provides some indications in this direction.

	\bigskip\bigskip

{\bf 1.1. The basic frame.} A pair of exact sequences 
$$
 0 \to U_0 @>w_0>> U_1 \to H \to 0 \t{and} 0 \to K\to U_0 @>v_0>> U_1 \to Q \to 0
$$
yields a module $U_2$ and a pair of exact sequences
$$
 0 \to U_1 @>w_1>> U_2 \to H \to 0 \t{and} 0 \to K \to U_1 @>v_1>> U_2 \to Q \to 0
$$
by forming the induced exact sequence of $ 0 \to U_0 @>w_0>> U_1 \to H \to 0 $
using the map $v_0$:
$$
\CD
   @.   0   @.   0 \cr
 @.     @VVV    @VVV         \cr
   @.   K   @=   K \cr
 @.     @VVV    @VVV         \cr
 0 @>>> U_0 @>w_0>> U_1 @>>> H @>>> 0 \cr
 @.     @VVv_0V    @VVv_1V   @|       \cr
 0 @>>> U_1 @>w_1>> U_2 @>>> H @>>> 0 \cr
 @.     @VVV    @VVV         \cr
   @.   Q   @=   Q \cr
 @.     @VVV    @VVV         \cr
   @.   0   @.   0 \cr
\endCD
$$

Recall that a commutative square
$$
\CD X @>f>> Y_1 \cr
  @VgVV @VVg'V \cr
  Y_2 @>>f'> Z
\endCD
$$
is said to be {\it exact} provided it is both a pushout and a pullback, thus if
and only if the sequence
$$
 0 \to X @>{\left[\smallmatrix f \cr g \endsmallmatrix \right]}>> Y_1\oplus Y_2
 @>{\left[\smallmatrix g' & -f' \endsmallmatrix \right]}>> Z \to 0
$$
is exact. Note that our basic setting provides an exact square
$$
\CD U_0 @>w_0>> U_1 \cr
  @Vv_0VV @VVv_1V \cr
  U_1 @>>w_1> U_2
\endCD
$$

Next, we will use that the composition of exact squares is exact:
\item{(E1)} {\it The composition of two exact squares 
$$
\CD X @>>> Y_1 @>>> Z_1 \cr
  @VVV @VVV   @VVV\cr
  Y_2 @>>> Z_2 @>>> A
\endCD
$$
yields an exact square}
$$
\CD X @>>>  Z_1 \cr
  @VVV   @VVV\cr
  Y_2 @>>> A
\endCD
$$

	\bigskip
{\bf 1.2. The ladder.} Using induction, we obtain in this way modules $U_i$ and pairs of exact sequences 
$$
 0 \to U_i @>w_i>> U_{i+1} \to H \to 0 \t{and} 0 \to K\to U_i @>v_i>> U_{i+1} \to Q \to 0
$$
for all $i \ge 0.$ 
	\medskip
We may combine the pushout diagrams constructed inductively and obtain the following 
ladder of commutative squares:
$$
\CD
 U_0 @>w_0>> U_1 @>w_1>> U_2 @>w_2>> U_3 @>w_3>> \cdots \cr
 @Vv_0VV   @Vv_1VV    @Vv_2VV    @Vv_3VV      \cr   
 U_1 @>w_1>> U_2 @>w_2>> U_3 @>w_3>>  U_4 @>w_4>> \cdots \cr
\endCD
$$
We form the inductive limit $U_\infty = \bigcup_i U_i$ 
(along the maps $w_i$).
	\medskip
Since all the squares commute, the maps $v_i$ induce a map  $U_\infty  \to 
U_\infty $
which we denote by $v_\infty$:
$$
\CD
 U_0 @>w_0>> U_1 @>w_1>> U_2 @>w_2>> U_3 @>w_3>> \cdots  \qquad  \bigcup_i U_i  = U_\infty \cr
 @Vv_0VV   @Vv_1VV    @Vv_2VV    @Vv_3VV     \qquad           @VVv_\infty V \cr   
 U_1 @>w_1>> U_2 @>w_2>> U_3 @>w_3>>  U_4 @>w_4>> \cdots \qquad  \bigcup_i U_i = U_\infty \cr
\endCD
$$
We also may consider the factor modules $U_\infty/U_0$ and  $U_\infty/U_1$.
The map $v_\infty\:U_\infty \to U_\infty$ maps $U_0$ into $U_1$, thus it induces
a map 
$$
 \overline {v}\:U_\infty/U_0 \longrightarrow U_\infty/U_1.
$$

Claim. {\it The map $\overline{v}$ is an isomorphism.} Namely, the commutative diagrams
$$
\CD
 0 @>>> U_{i-1} @>w_{i-1}>> U_{i} @>>> H @>>> 0 \cr
 @.     @VVv_{i-1}V    @VVv_{i}V   @|       \cr
 0 @>>> U_{i} @>w_{i}>> U_{i+1} @>>> H @>>> 0 \cr
\endCD
$$
can be rewritten as 
$$
\CD
 0 @>>> U_{i-1} @>w_{i-1}>> U_{i} @>>> U_{i}/U_{i-1} @>>> 0 \cr
 @.     @VVv_{i-1}V    @VVv_{i}V   @VV{\overline{v}_i}V       \cr
 0 @>>> U_{i} @>w_{i}>> U_{i+1} @>>> U_{i+1}/U_{i} @>>> 0 \cr
\endCD
$$
with an isomorphism $\overline{v}_i\: U_{i}/U_{i-1} \to U_{i+1}/U_{i}.$
The map
$\overline{v}$ is a map from a filtered module with factors $U_i/U_{i-1}$ (where $i\ge 1$)
to a filtered module with factors $U_{i+1}/U_i$ (again with $i\ge 1$), and the maps
$\overline{v_i}$ are just those induced on the factors. 

It follows: The composition of maps
$$
\CD
 U_\infty/U_0  @>p>> U_\infty/U_1   @>{\overline {v}^{-1}}>>   U_\infty/U_0
\endCD
$$
with $p$ the projection map is an epimorphism $\phi$ with kernel $U_1/U_0.$
It is easy to see that $\phi$ is locally nilpotent, namely we have $\phi^t(U_t/U_0) = 0$
for all $t$.
	\bigskip
{\bf Summery. (a)} {\it The maps $v_i$ yield a map
$$
  v_\infty\: U_\infty \to U_\infty
$$ 
with kernel $K$ and cokernel $Q$.

{\bf (b)}
This map $v_\infty$ induces an isomorphism $\overline{v}\: U_\infty/U_0 \to U_\infty/U_1$.
Composing the inverse of this isomorphism with the canonical projection $p$, we obtain
an endomorphism $\phi = (\overline v)^{-1}\circ p$
$$
 U_\infty/U_0 @>p>> U_\infty/U_1 @>\overline v^{-1}>> U_\infty/U_0.
$$
If the cokernel $H$ of $w$ 
is non-zero and of finite length, then 
$U_\infty/U_0$ 
is a Pr\"ufer module with respect to $\phi,$ with basis $H$;} in this case,
we call $U_\infty/U_0$ 
(or better the pair $(U_\infty/U_0,\phi)$)  the
{\it Pr\"ufer module defined by the pair $(w_0,v_0)$} or by the ladder $U_i$.
Pr\"ufer modules which are obtained in this way will be said to be {\it of ladder type.}
	\medskip
If necessary, we will use the following 
notation: $U_i(w_0,v_0) = U_i$, for all $i\in \Bbb N\cup \{\infty\}$
and 
$P(w_0,v_0) = U_\infty/U_0$ for the Pr\"ufer module. Since
$P(w_0,v_0)$ is a Pr\"ufer module with basis the cokernel $H$ of $w$, we will sometimes
write $H[n] = U_n/U_0$ or even $H[n;w_0,v_0].$ 
	\bigskip

{\bf Remark:} Using a terminology introduced for string algebras [R3], we also could say:
$U_\infty$ is {\it expanding,} $U_\infty/U_0$ is {\it contracting.}
	\bigskip
{\bf Lemma.} 
{\it Assume that $P = P(w,v)$ with $w,v\:U_0 \to U_1$. Then $P$ is generated by $U_1$.}
	\medskip
thus by induction $U_i$ is a factor module of the direct sum of $i$ copies of $U_1.$  

	\bigskip
%=================================================================
{\bf 1.3. The chessboard.} Assume now that both maps $w_0,v_0\: U_0 \to U_1$ are
monomorphisms. Then we get the following arrangement of commutative squares:
$$
\CD
 U_0 @>w_0>> U_1 @>w_1>> U_2 @>w_2>> U_3 @>w_3>> \cdots \cr
 @Vv_0VV   @Vv_1VV    @Vv_2VV    @Vv_3VV      \cr   
 U_1 @>w_1>> U_2 @>w_2>> U_3 @>w_3>> \cdots    \cr
 @Vv_1VV   @Vv_2VV    @Vv_3VV         \cr   
 U_2 @>w_2>> U_3 @>w_3>>  \cdots    \cr
 @Vv_2VV   @Vv_3VV              \cr   
 U_3 @>w_3>> \cdots    \cr
 @Vv_3VV                 \cr   
 \cdots    \cr
\endCD
$$
Note that there are 
both horizontally as well as vertically ladders: the horizontal ladders yield
$U_\infty(w_0,v_0)$ (and its endomorphism $v_\infty$); the vertical ladders yield
$U_\infty(v_0,w_0)$ (and its endomorphism $w_\infty$). 
	\bigskip\bigskip
%==============================================================================
{\bf 2. Examples.}
	\medskip
{\bf (1)} The classical example: Let $R = \Bbb Z$ be the ring of integers, 
and $U_0 = U_1 = \Bbb Z$  its regular representation. Module homomorphisms
$\Bbb Z \to \Bbb Z$ are given by the multiplication with some integer $n$, thus we denote
such a map just by $n$. Let $w_0 = 2$ and $v_0 = n$. {\it If $n$ is odd, then $P(2,n)$ 
is the ordinary Pr\"ufer group for the prime $2$,} and $U_\infty(2,n) = \Bbb Z[\frac12]$
(the subring of $\Bbb Q$ generated by $\frac12$). 
If $n$ is even, then $P(2,n)$ is an elementary abelian 2-group. 
	\medskip
{\bf (2)} Let $R = K(2)$ be the Kronecker algebra over some field $k$.
Let $U_0$ be simple projective, $U_1$ indecomposable projective of length 3 and
$w_0\:U_0\to U_1$ a non-zero map with cokernel $H$ (one of the indecomposable modules of length
2). The module $P(w_0,v_0)$ is the Pr\"ufer module for $H$ if and only if $v_0 \notin
kw_0,$ otherwise it is a direct sum of copies of $H$.
	\medskip
{\bf (3)} Trivial cases: First, let $w$ be a split monomorphism. Then the Pr\"ufer module
with respect to any map $v\:U_0 \to U_1$ is just the countable sum
of copies of $H$. Second, let $w\:U_0 \to U_1$ be an arbitrary monomorphism, let
$\beta\:U_1 \to U_1$ be an endomorphism. Then $P(w,\beta w)$ is the countable sum
of copies of $H$.
	\medskip
{\bf (4)} Assume that there exists a split monomorphism $v\: U_0 \to U_1$,
say $U_1 = U_0 \oplus X$ and 
$v = \left[\smallmatrix 1 \cr 0\endsmallmatrix\right]\:U_0 \to U_1$.
Then
$$
 0 \to U_0 @>w>> U_0\oplus X @>>> H \to 0
$$
is a Riedtmann-Zwara sequence as discussed in section 4, 
thus $H$ is a degeneration of $X$.
	\bigskip
{\bf Remark:} {\it Not all Pr\"ufer modules are of ladder type.} Consider the generalized
Kronecker algebra $\Lambda$ with 
countably many arrows $\alpha_0,\alpha_1,\dots$ starting at the
vertex $a$ and ending in the vertex $b$. 
Define a representation $P = (P_a,P_b,\alpha_i)_i$ 
as follows: Let $P_a = P_b$ be a vector space with a countable basis
$e_0,e_1,\dots$ and let $\alpha_i\:P_a \to P_b$ be defined by $\alpha_i(e_j) = e_{j-i}$
provided $j\ge i$ and $\alpha_i(e_j)= 0$ otherwise. Let $\phi_a, \phi_b$ 
be the endomorphism of $P_a$ of $P_b$, respectively, which sends $e_0$ to $0$ and $e_i$ to 
$e_{i-1}$ for $i \ge 1.$ Then $P$ is a Pr\"ufer module (with respect to $\phi$, but also with
respect to any power of $\phi$). 
Obviously, $P$ is a faithful
$\Lambda$-module. Assume that $P = P(w,v)$ for some maps $w,v\: U_0 \to U_1$ with $U_0,U_1$
of finite length. Then $P$ is
generated by $U_1$, according to Lemma 1.2. However $U_1$ is of finite length and no 
finite length $\Lambda$-module is faithful.
	\bigskip\bigskip
%============================================================================
{\bf 3. Ladder extensions.}
	\medskip
{\bf 3.1. The definition.}
A self-extension $ 0 \to H \to H[2] \to H \to 0$ is said to be a {\it ladder extension} provided
there is a commutative diagram with exact rows 
$$
\CD
 0 @>>> U_0 @>>> U_1 @>q>> H @>>> 0 \cr
 @.     @VVf V    @VVV   @|       \cr
 0 @>>> H @>>> H[2] @>>> H @>>> 0 \cr
\endCD
$$
such that $f$ factors through $q$, say $f = qv$ for some $v\:U_0 \to U_1.$ 
In case $U_0$ is in addition a simple module, we say that $\epsilon$ is of {\it simple}
ladder type.
	\bigskip
This means that we have a commutative diagram with exact rows of the following kind
(here $f = qv_0$):
$$
\CD 
 0 @>>> U_0 @>w_0>> U_1 @>q>> H @>>> 0 \cr
 @.    @Vv_0VV    @Vv_1VV     @| \cr
 0 @>>> U_1 @>{w_1}>> U_2 @>>> H @>>> 0 \cr
 @.    @VqVV    @VVV     @| \cr
 0 @>>> H @>>> H[2] @>>> H @>>> 0.
\endCD
$$
Thus, in order to construct all the ladder extensions of $H$, we may start with
an arbitrary epimorphism $q\:U_1 \to H,$ form its kernel $w_0$ and consider
any homomorphism $v_0\: U_0\to U_1$.
	\medskip
	
According to section 1 we know: {\it Ladder extensions built up to form Pr\"ufer
modules.}
	\bigskip

{\bf Lemma.} {\it Let $k$ be a commutative ring and $\Lambda$ a $k$-algebra.
Then $H[2;w_0,v_0] = H[2;w_0,v_0+\mu w_0]$ for any $\mu \in k.$}
	\medskip
Proof: We deal with the exact sequence induced by $qv_0$ or $q(v_0+\mu w_0)$,
respectively. But $q(v_0+\mu w_0) = qv_0 + q\mu w_0 = qv_0,$ since $qw_0 = 0.$
	\medskip
Also, any central automorphism $\lambda$ of $U_0$ yields isomorphic
extensions $H[2;w_0,v_0]$ and $H[2;w_0\lambda v_0]$. This shows that {\it the extension
$H[2;w_0,v_0]$ only depends on the $k$-subspace $\langle w_0, v_0\rangle.$}
	\bigskip
{\bf Remark.} Not all self-extensions are ladder extensions. For example:
{\it 
A non-zero self-extension of a simple module $S$ over an artinian ring is never a ladder
extension!}
	\medskip
Proof: Construct the corresponding ladder, thus the corresponding Pr\"ufer module $S[\infty]$.
The module $S[n]$ would be a (serial) module of Loewy length $n$, with $n$ arbitrary. But the
Loewy length of any module over the artinian ring $R$ is bounded by the Loewy length of
${}_RR,$ thus $S[\infty]$ cannot exist.
	\medskip
%===The Auslander question=========
	\medskip
{\bf Example.} Here is a further example of a self-extension which is not a ladder extension.
Consider the following quiver $Q$
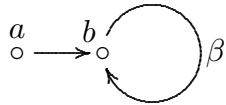
$$
{\beginpicture
\setcoordinatesystem units <1cm,1cm>
\put{} at 0 0.5 
\put{} at 1.5 -0.5 
\put{$a$} at 0 0.25 
\put{$b$} at 0.8 0.25 
\put{$\beta$} at 2.2 0 
\put{$\circ$} at 0 0
\put{$\circ$} at 0.95 0
\arr{0.2 0}{0.8 0}
\circulararc -320 degrees from 1 0.2  center at 1.5 0 
\arr{1.01 -0.2}{1 -0.18}
\endpicture}
$$
with one loop $\beta$ at the vertex b, and one arrow from $a$ to $b$. We consider the
representations of $Q$ with the relation $\beta^3 = 0.$
The universal covering $\widetilde Q$ of $Q$ has many $D_5$ subquivers $Q'$ of the form
$$
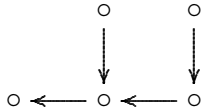
{\beginpicture
\setcoordinatesystem units <1cm,1cm>
\multiput{$\circ$} at  0 0  1 0  1 1  2 0  2 1 /
\arr{1 0.8}{1 0.2}
\arr{2 0.8}{2 0.2}

\arr{0.8 0}{0.2 0}
\arr{1.8 0}{1.2 0}
\endpicture} 
$$
and we consider some representations of $Q'$; we present here the corresponding dimension
vectors. 
$$
{\beginpicture
\setcoordinatesystem units <1cm,1cm>
\put{\beginpicture
\put{$0$} at 0 1
\multiput{$1$} at 0 0  1 1  2 0  2 1 /
\put{$2$} at 1 0
\arr{0 0.8}{0 0.2}
\arr{1 0.8}{1 0.2}
\arr{2 0.8}{2 0.2}

\arr{0.8 0}{0.2 0}
\arr{1.8 0}{1.2 0}
\put{$H'$} at 1 -0.5 
\endpicture} at 0 0 

\put{\beginpicture
\multiput{$0$} at 0 1  2 0  2 1 /
\multiput{$1$} at 0 0  1 0  1 1 /
\arr{0 0.8}{0 0.2}
\arr{1 0.8}{1 0.2}
\arr{2 0.8}{2 0.2}

\arr{0.8 0}{0.2 0}
\arr{1.8 0}{1.2 0}
\put{$H$} at 1 -0.5 
\endpicture} at -4 0 

\put{\beginpicture
\multiput{$0$} at 0 0  0 1  1 1 /
\multiput{$1$} at 1 0  2 0  2 1 /
\arr{0 0.8}{0 0.2}
\arr{1 0.8}{1 0.2}
\arr{2 0.8}{2 0.2}

\arr{0.8 0}{0.2 0}
\arr{1.8 0}{1.2 0}
\put{$H''$} at 1 -0.5 
\endpicture} at 4 0 

\endpicture}
$$
There is an obvious exact sequence
$$
 0 \to H \to H' \to H'' \to 0.
$$
Under the covering functor, the representations $H$ and $H''$ are identified, thus
we obtain a self-extension.
One easily checks that this self-extension is not a ladder extension. 
	\bigskip
{\bf Proposition.} {\it Let $H$ be an indecomposable module with Auslander-Reiten translate isomorphic to $H$. Assume that there is a simple submodule $S$ of $H$ with $\Ext^1(S,S) = 0.$ Then the
Auslander-Reiten sequence ending (and starting) in $H$ is a ladder extension.}
	\medskip
Proof. Let $0 @>>> H @>>> H' @>>> H @>>> 0$ be the Auslander sequence.
Denote by $u\: S \to H$ the inclusion map. Since the maps $H \to H/S$ factors through 
$H \to H'$,
there is a commutative diagram with exact rows of the following form:
$$
\CD 
 0 @>>>  S @>w>> U @>q>> H @>>> 0 \cr
 @.    @Vu VV    @VVV     @| \cr
 0 @>>>   H     @>>>  H' @>>> H @>>> 0
\endCD
$$
Now form he induced exact sequence:
$$
\CD 
 0 @>>>   S    @>>>  U' @>>> S @>>> 0 \cr
 @.    @|    @VVV     @Vu VV \cr
 0 @>>>  S @>w>> U @>q>> H @>>> 0 
\endCD
$$
Since $\Ext^1(S,S) = 0,$ the induced sequence splits, thus 
we obtain a map $v\:S \to U$ with $qv = u.$ It follows that $H' = H[2;w,v].$

	\medskip
We do not know whether one can delete the assumption about the existence of $S$. 
	\bigskip
%===============================================================================
{\bf 3.2. Standard self-extensions.}
	\medskip
Let $H$ be an $R$-module, say with an exact sequence 
$0 \to \Omega H @>u>> PH @>p>> H \to 0$, where $PH$ denotes a projective cover of $H$.
We know that 
$$
 \Ext^1(H,H) = \underline{\Hom}(\Omega H,H) = \Hom(\Omega H,H)/\Im(\Hom(u,H)).
$$
Note that 
$$
 \Im(\Hom(u,H)) \subseteq \Im(\Hom(\Omega H,p)) \subseteq \Hom(\Omega H,H).
$$
(Proof: $\Hom(u,H)\:\Hom(PH,H) \to \Hom(\Omega H,H)$, thus take $\phi\:PH \to H$
 and form $\phi u.$ Since $p\:PH \to H$ is surjective and $PH$ is projective,
 there is $\phi'\:PH \to PH$ with $\phi = p\phi'.$
 Thus $\phi u = p\phi' u$ is in the image of $\Hom(\Omega H,p).$) 
	
Thus we can consider
$$
 \Ext^1(H,H)_s := \Im(\Hom(\Omega H,p))/\Im(\Hom(u,H))
$$
as a subgroup of $\Hom(\Omega H,H)/\Im(\Hom(u,H)) =
 \Ext^1(H,H).$ We call the elements of $\Ext^1(H,H)_s$ the {\it standard } self-extensions.
	\bigskip
{\bf Proposition.} {\it Standard self-extensions are ladder extensions.}
	\medskip
{\bf Proof.}
Here is the usual diagram in which way a map $f\:\Omega H \to H$ yields an
self-extension of $H$
$$
\CD 
 0 @>>> \Omega H @>u>> PH @>p>> H @>>> 0 \cr
 @.    @Vf VV    @VVV     @| \cr
 0 @>>>   H      @>>>  H[2] @>>> H @>>> 0
\endCD
$$
The standard  extensions are those where the map $f$ factors through $p$, say
$f = pw'$ with $w'\:\Omega H \to PH:$
$$
\CD 
 0 @>>> \Omega H @>u>> PH @>p>> H @>>> 0 \cr
   @.    @Vw' VV    @Vw'_1VV     @| \cr
0 @>>>  PH      @>u_1>>  U_2 @>>> H @>>> 0 \cr
 @.    @Vp VV    @VVV     @| \cr
 0 @>>>   H      @>>>  H[2] @>>> H @>>> 0
\endCD
$$
	\bigskip
%===============================================================================
{\bf 3.3. Modules of projective dimension 1.}
	\medskip
{\bf Proposition.} {\it If the projective dimension of $H$ is at most 1, then any
self-extension of $H$ is standard, thus a ladder extension.}
	 \medskip
Proof: Consider a module $H$ with a projective presentation $0 \to P' \to P
@>p>> H \to 0.$ Any self-extension of $H$ is given by a diagram of the
following kind:
$$
\CD 
 0 @>>> P' @>u>> P @>p>> H @>>> 0 \cr
 @.    @Vf VV    @VVV     @| \cr
 0 @>>>   H      @>>>  H[2] @>>> H @>>> 0
\endCD
$$
Since $P'$ is projective and $p\:P \to H$ surjective, there is a map $
f'\:P' \to P$ such that $f = pf'.$ The self-extension is 
given just by $H[2] = H[2;u,f'].$
	\medskip
{\bf Corollary.} 
{\it If $R$ is a hereditary ring, any self-extension is standard, thus a ladder extension.}
	\bigskip
{\bf Example} of a ladder extension which is not standard.
Consider the quiver $Q$
$$
{\beginpicture
\setcoordinatesystem units <1cm,1cm>
\put{} at  0.2 0
\put{} at -0.2 2
\put{$c$} at 0.25 -0.1
\put{$b$} at  0.3 1 
\put{$a$} at 0.25 2.1
\put{$\circ$} at  0 0 
\put{$\circ$} at 0 1
\put{$\circ$} at 0 2
\arr{0.1 0.8}{0.1 0.2}
\arr{-0.1 0.8}{-0.1 0.2}
\arr{0.1 1.8}{0.1 1.2}
\arr{-0.1 1.8}{-0.1 1.2} 
\put{$\alpha$} at -0.3 1.5 
\put{$\beta$} at  0.3 1.5
\put{$\gamma$} at -0.3 0.5
\put{$\delta$} at 0.3 0.5
\endpicture}
$$
such that $\delta\alpha = 0 = \gamma\beta = \gamma\alpha-\delta\beta.$ 
Consider the indecomposable length 2 module $H =  (\beta\:a \to b)$ annihilated by $\alpha$.
Then the kernel $\Omega H$ of $PH \to H$ is $\Omega H = (\gamma\:b \to c).$ We may visualize
this as follows:
$$
{\beginpicture
\setcoordinatesystem units <.5cm,.5cm>
\multiput{$\bullet$} at 0 1  1 0  6 1  7 0  7 2  8 1   13 2  14 1 /
\plot 0 1  1 0  /
\plot 7 0  6 1  7 2  8 1  7 0 /
\plot 13 2  14 1 /
\arr{2 1}{4.5 1}
\arr{9.5 1}{12 1}
\arr{-3 1}{-1 1}
\arr{15 1}{17 1}
\put{$0$} at -4 1
\put{$0$} at 19 1
\put{$\gamma$} at 0.2 0.2 
\put{$\gamma$} at 6.2 0.2 
\put{$\delta$} at 7.8 0.2 
\put{$\beta$} at 7.8 1.8 
\put{$\alpha$} at 6.2 1.8 
\put{$\beta$} at 13.8 1.8 
\put{$\sssize b$} at 0 1.4
\put{$\sssize c$} at 1 -.4
\put{$\sssize a$} at 7 2.4
\put{$\sssize b$} at 5.7 1
\put{$\sssize b$} at 8.3 1
\put{$\sssize c$} at 7 -0.4
\put{$\sssize a$} at 13 2.4
\put{$\sssize b$} at 14 0.6
\put{$u$} at 3 1.5
\put{$p$} at 10.5 1.5
\endpicture}
$$
There is a ladder extension of $H$, given by the non-trivial map 
$f\:\Omega H \to H$,
but this map does not factor through $PH$, since $\Hom(\Omega H,PH)$ is one-dimensional, generated by $u$. Note that $\overline f\:\Omega H/K$ factors through 
$\overline p\:PH/u(K) \to H$, where $K = S(c)$ is the kernel of $f$. 
	\bigskip

%=========================================================================
{\bf Appendix.} Here, we want to indicate hat Corollary can be used in order to obtain
a conceptual proof of the second Brauer-Thrall conjecture for hereditary artin algebras.
	\medskip
Assume that there is no generic module. We show: {\it Any indecomposable module is a brick without 
self extensions.}
Assume that there is an indecomposable module $M$ which is not a brick of which does have 
self-extensions. If $M$ is not a brick, then the brick paper [R2] shows that there
are bricks $M'$ with self-extensions. Thus, we see that there always is a brick 
$H$ with self-extensions. Take any non-zero self-extension of $H$. 
According to 3.2, such a self-extension is standard, thus a ladder
extension,
thus we obtain a corresponding Pr\"ufer module $H[\infty]$. The process of simplification 
[R1] shows that
all the modules $H[n]$ are indecomposable. Thus $H[\infty]$ is not of finite type and therefore
there exists a generic module [R5].

But if any indecomposable module is a brick without self-extensions, the quadratic form is
weakly positive. Ovsienko asserts that then there are only finitely many positive roots,
thus the algebra is of bounded representation type and therefore of finite representation type. 
	\bigskip\bigskip
%================================================================================
{\bf 3.4. Warning.} 
{\it A Pr\"ufer module $M[\infty]$ is not necessarily determined by $M[2]$, even if it
is of ladder type.}
	\medskip
As an example take the generalized Kronecker quiver with vertices $a,b$ and three arrows 
$\alpha,\beta,\gamma\:a \to b$.
and let $H$ be the two-dimensional indecomposable representation annihilated by $\alpha$
and $\beta$. Consider a projective cover $q\:PH \to H$, let $\Omega H$ be its kernel,
say with inclusion map $w\:\Omega H \to PH.$
$$
 0 @>>> \Omega H @>w>> PH @>q>> H @>>> 0 \tag{*}
$$
The ladders to be considered are given by the various maps $f\:\Omega H \to PH$ such that
the image of $f$ is not contained in $\Omega H$ (otherwise, the induced self-extension of $H$ will split). In order to specify a self-extension $H[2]$ of $H$, we require that $H[2]$ is
annihilated say by $\gamma$. 

We will consider several copies of $PH$. If $e_i \in (PH)_a$ is a generator, let us
denote $e_{i1} = \alpha(e_i),\ 
e_{i2} = \beta(e_i),\ e_{i3} = \gamma(e_i),$ thus, $e_{i1}, e_{i2}, e_{i3}$ is a basis
of $(PH)_b$. 

We start with $PH$ generated by $e_1$ and consider the exact sequence $(*)$ as displayed 
above. We see that $e_{12}, e_{13}$ is a basis of $\Omega H.$
 
Now, let us consider two maps $f,g\:\Omega H \to PH$, here we denote the generator of $PH$
by $e_0.$ The first map $f$ is given
by $f(e_{12}) = e_{01}$ and $f(e_{13}) = 0.$
The second map $g\:\Omega H \to PH$ is defined by $g(e_{12}) = e_{01}$ and 
$g(e_{13}) = e_{02}$. 

Note that $qf = qg$, thus $H[2;w,f] = H[2;w,g]$ and actually this is precisely the self-extension
of $H$ annihilated by $\gamma.$

An easy calculation shows that $H[3;w,f]$ (and even $H[\infty;w,f]$) is annihilated by $\gamma$, whereas $H[3;w,g]$ us faithful. The following displays may be helpful; always, we exhibit
the modules:
$$
\CD
U_0 = \Omega H @>w_0>> U_1 = PH @>w_1>> U_2 \cr
   @VVv_0V    @VVv_1V    @VVv_2V   \cr     
U_1 = PH @>w_1>> U_2 @>w_2>> U_3 \cr
   @VVq V    @VVV    @VVV   \cr     
  H @>>> H[2] @>>> H[3] \cr
\endCD
$$
First the display for the homomorphism $f$.
$$
{\beginpicture
\setcoordinatesystem units <.6cm,.6cm>
%============\Omega H==============================
\put{\beginpicture
\put{} at 1 1
\put{$\ssize e_{13}$} at 0 0
\put{$\ssize e_{12}$} at 1 0
\endpicture} at 0 8
%==========PH==============================
\put{\beginpicture
\put{$\ssize e_1$} at 1 1
\put{$\ssize e_{13}$} at 0 0
\put{$\ssize e_{12}$} at 1 0
\put{$\ssize e_{11}$} at 2 0
\plot 0.8 0.8  0.2 0.2 /
\plot 1 0.8  1 0.2 /
\plot 1.2 0.8  1.8 0.2 /
\endpicture} at 8 8
%============U_2==============================
\put{\beginpicture
\put{$\ssize e_1$} at 1 1
\put{$\ssize e_{13}$} at 0 0
\put{$\ssize e_{12}$} at 1 0
\put{$\ssize e_{11}$} at 2 0

\put{$\ssize e_2$} at 2 1
\put{$\ssize e_{22}$} at 2 -0.3
\put{$\ssize e_{21}$} at 3 0
\plot 0.8 0.8  0.2 0.2 /
\plot 1 0.8  1 0.2 /
\plot 1.2 0.8  1.8 0.2 /
\plot 2 0.8  2 0.2 /
\plot 2.2 0.8  2.8 0.2 /

\endpicture} at 16 8

%============PH==============================
\put{\beginpicture
\put{$\ssize e_0$} at 1 1
\put{$\ssize e_{03}$} at 0 0
\put{$\ssize e_{02}$} at 1 0
\put{$\ssize e_{01}$} at 2 0
\plot 0.8 0.8  0.2 0.2 /
\plot 1 0.8  1 0.2 /
\plot 1.2 0.8  1.8 0.2 /
\endpicture} at 0 4

%============U_2 unten==============================
\put{\beginpicture
\put{$\ssize e_0$} at 1 1
\put{$\ssize e_{03}$} at 0 0
\put{$\ssize e_{02}$} at 1 0
\put{$\ssize e_{01}$} at 2 0

\put{$\ssize e_1$} at 2 1
\put{$\ssize e_{12}$} at 2 -0.3
\put{$\ssize e_{11}$} at 3 0
\plot 0.8 0.8  0.2 0.2 /
\plot 1 0.8  1 0.2 /
\plot 1.2 0.8  1.8 0.2 /
\plot 2 0.8  2 0.2 /
\plot 2.2 0.8  2.8 0.2 /

\endpicture} at 8 4

%============U_3 unten==============================
\put{\beginpicture
\put{$\ssize e_0$} at 1 1
\put{$\ssize e_{03}$} at 0 0
\put{$\ssize e_{02}$} at 1 0
\put{$\ssize e_{01}$} at 2 0

\put{$\ssize e_1$} at 2 1
\put{$\ssize e_{12}$} at 2 -0.3
\put{$\ssize e_{11}$} at 3 0

\put{$\ssize e_2$} at 3 1
\put{$\ssize e_{22}$} at 3 -0.3
\put{$\ssize e_{21}$} at 4 0

\plot 0.8 0.8  0.2 0.2 /
\plot 1 0.8  1 0.2 /
\plot 1.2 0.8  1.8 0.2 /
\plot 2 0.8  2 0.2 /
\plot 2.2 0.8  2.8 0.2 /
\plot 3 0.8  3 0.2 /
\plot 3.2 0.8  3.8 0.2 /

\endpicture} at 16 4
%===========
%============H==============================
\put{\beginpicture
\put{$\ssize e_0$} at 1 1
\put{} at 0 0

\put{$\ssize e_{01}$} at 2 0
\plot 1.2 0.8  1.8 0.2 /
\endpicture} at 0 0

%============H[2] unten==============================
\put{\beginpicture
\put{$\ssize e_0$} at 1 1
\put{} at 0 0
\put{$\ssize e_{01}$} at 2 0

\put{$\ssize e_1$} at 2 1
\put{$\ssize e_{12}$} at 2 -0.3
\put{$\ssize e_{11}$} at 3 0
\plot 1.2 0.8  1.8 0.2 /
\plot 2 0.8  2 0.2 /
\plot 2.2 0.8  2.8 0.2 /

\endpicture} at 8 0

%============H[3] unten==============================
\put{\beginpicture
\put{$\ssize e_0$} at 1 1
\put{} at 0 0
\put{$\ssize e_{01}$} at 2 0

\put{$\ssize e_1$} at 2 1
\put{$\ssize e_{12}$} at 2 -0.3
\put{$\ssize e_{11}$} at 3 0

\put{$\ssize e_2$} at 3 1
\put{$\ssize e_{22}$} at 3 -0.3
\put{$\ssize e_{21}$} at 4 0

\plot 1.2 0.8  1.8 0.2 /
\plot 2 0.8  2 0.2 /
\plot 2.2 0.8  2.8 0.2 /
\plot 3 0.8  3 0.2 /
\plot 3.2 0.8  3.8 0.2 /

\endpicture} at 16 0
\arr{3 8}{5 8}
\arr{11 8}{13 8}
\arr{3 4}{5 4}
\arr{11 4}{13 4}
\arr{3 0}{5 0}
\arr{11 0}{13 0}

\arr{0 6.5}{0 5.5}
\arr{8 6.5}{8 5.5}
\arr{16 6.5}{16 5.5}

\arr{0 2.5}{0 1.5}
\arr{8 2.5}{8 1.5}
\arr{16 2.5}{16 1.5}
\endpicture}
$$

\noindent
Now  the corresponding display for the homomorphism $g$.
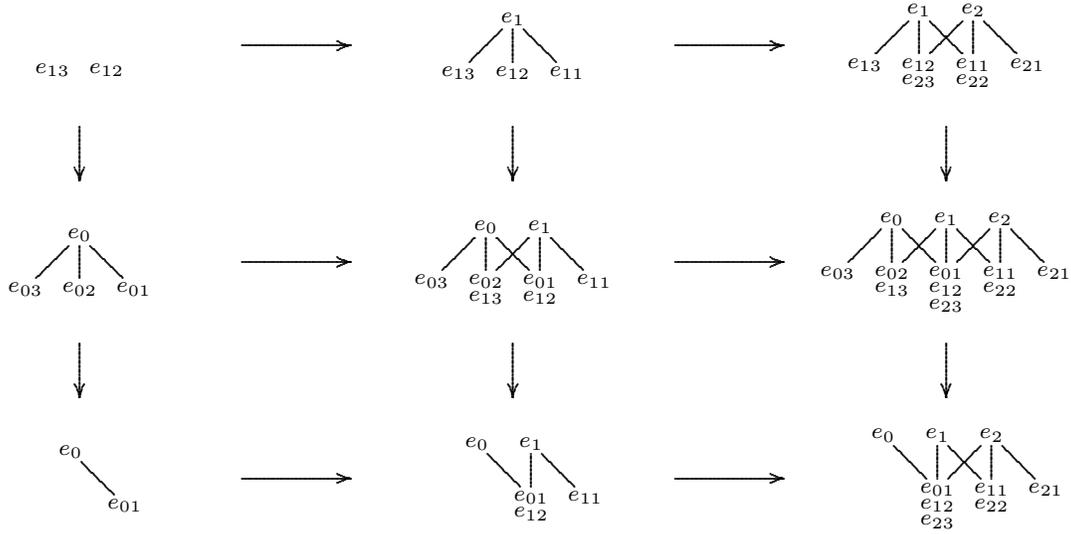
$$
{\beginpicture
\setcoordinatesystem units <.6cm,.6cm>
%============\Omega H==============================
\put{\beginpicture
\put{} at 1 1
\put{$\ssize e_{13}$} at 0 0
\put{$\ssize e_{12}$} at 1 0
\endpicture} at 0 8
%==========PH==============================
\put{\beginpicture
\put{$\ssize e_1$} at 1 1
\put{$\ssize e_{13}$} at 0 0
\put{$\ssize e_{12}$} at 1 0
\put{$\ssize e_{11}$} at 2 0
\plot 0.8 0.8  0.2 0.2 /
\plot 1 0.8  1 0.2 /
\plot 1.2 0.8  1.8 0.2 /
\endpicture} at 8 8
%============U_2==============================
\put{\beginpicture
\put{$\ssize e_1$} at 1 1
\put{$\ssize e_{13}$} at 0 0
\put{$\ssize e_{12}$} at 1 0
\put{$\ssize e_{11}$} at 2 0

\put{$\ssize e_2$} at 2 1
\put{$\ssize e_{23}$} at 1 -0.3
\put{$\ssize e_{22}$} at 2 -0.3
\put{$\ssize e_{21}$} at 3 0
\plot 0.8 0.8  0.2 0.2 /
\plot 1 0.8  1 0.2 /
\plot 1.2 0.8  1.8 0.2 /
\plot 2 0.8  2 0.2 /
\plot 2.2 0.8  2.8 0.2 /
\plot 1.8 0.8  1.2 0.2 /
\endpicture} at 16 8

%============PH==============================
\put{\beginpicture
\put{$\ssize e_0$} at 1 1
\put{$\ssize e_{03}$} at 0 0
\put{$\ssize e_{02}$} at 1 0
\put{$\ssize e_{01}$} at 2 0
\plot 0.8 0.8  0.2 0.2 /
\plot 1 0.8  1 0.2 /
\plot 1.2 0.8  1.8 0.2 /
\endpicture} at 0 4

%============U_2 unten==============================
\put{\beginpicture
\put{$\ssize e_0$} at 1 1
\put{$\ssize e_{03}$} at 0 0
\put{$\ssize e_{02}$} at 1 0
\put{$\ssize e_{01}$} at 2 0

\put{$\ssize e_1$} at 2 1
\put{$\ssize e_{13}$} at 1 -0.3
\put{$\ssize e_{12}$} at 2 -0.3
\put{$\ssize e_{11}$} at 3 0
\plot 0.8 0.8  0.2 0.2 /
\plot 1 0.8  1 0.2 /
\plot 1.2 0.8  1.8 0.2 /
\plot 2 0.8  2 0.2 /
\plot 2.2 0.8  2.8 0.2 /
\plot 1.8 0.8  1.2 0.2 /

\endpicture} at 8 4

%============U_3 unten==============================
\put{\beginpicture
\put{$\ssize e_0$} at 1 1
\put{$\ssize e_{03}$} at 0 0
\put{$\ssize e_{02}$} at 1 0
\put{$\ssize e_{01}$} at 2 0

\put{$\ssize e_1$} at 2 1
\put{$\ssize e_{13}$} at 1 -0.3
\put{$\ssize e_{12}$} at 2 -0.3
\put{$\ssize e_{11}$} at 3 0

\put{$\ssize e_2$} at 3 1
\put{$\ssize e_{23}$} at 2 -0.6
\put{$\ssize e_{22}$} at 3 -0.3
\put{$\ssize e_{21}$} at 4 0

\plot 0.8 0.8  0.2 0.2 /
\plot 1 0.8  1 0.2 /
\plot 1.2 0.8  1.8 0.2 /
\plot 2 0.8  2 0.2 /
\plot 2.2 0.8  2.8 0.2 /
\plot 3 0.8  3 0.2 /
\plot 3.2 0.8  3.8 0.2 /
\plot 1.8 0.8  1.2 0.2 /
\plot 2.8 0.8  2.2 0.2 /
\endpicture} at 16 4
%===========
%============H==============================
\put{\beginpicture
\put{$\ssize e_0$} at 1 1
\put{} at 0 0

\put{$\ssize e_{01}$} at 2 0
\plot 1.2 0.8  1.8 0.2 /
\endpicture} at 0 0

%============H[2] unten==============================
\put{\beginpicture
\put{$\ssize e_0$} at 1 1
\put{} at 0 0
\put{$\ssize e_{01}$} at 2 0

\put{$\ssize e_1$} at 2 1
\put{$\ssize e_{12}$} at 2 -0.3
\put{$\ssize e_{11}$} at 3 0
\plot 1.2 0.8  1.8 0.2 /
\plot 2 0.8  2 0.2 /
\plot 2.2 0.8  2.8 0.2 /

\endpicture} at 8 0

%============H[3] unten==============================
\put{\beginpicture
\put{$\ssize e_0$} at 1 1
\put{} at 0 0
\put{$\ssize e_{01}$} at 2 0

\put{$\ssize e_1$} at 2 1
\put{$\ssize e_{12}$} at 2 -0.3
\put{$\ssize e_{11}$} at 3 0

\put{$\ssize e_2$} at 3 1
\put{$\ssize e_{23}$} at 2 -0.6
\put{$\ssize e_{22}$} at 3 -0.3
\put{$\ssize e_{21}$} at 4 0

\plot 1.2 0.8  1.8 0.2 /
\plot 2 0.8  2 0.2 /
\plot 2.2 0.8  2.8 0.2 /
\plot 3 0.8  3 0.2 /
\plot 3.2 0.8  3.8 0.2 /
\plot 2.8 0.8  2.2 0.2 /

\endpicture} at 16 0
\arr{3 8}{5 8}
\arr{11 8}{13 8}
\arr{3 4}{5 4}
\arr{11 4}{13 4}
\arr{3 0}{5 0}
\arr{11 0}{13 0}

\arr{0 6.5}{0 5.5}
\arr{8 6.5}{8 5.5}
\arr{16 6.5}{16 5.5}

\arr{0 2.5}{0 1.5}
\arr{8 2.5}{8 1.5}
\arr{16 2.5}{16 1.5}
\endpicture}
$$

	\bigskip\bigskip
%====================================
{\bf 4. Degenerations.}
	\medskip
Definition: Let $X, Y$ be finite length modules.
Call $Y$ a {\it degeneration} of $X$ provided there is an exact sequence
of the form $0 \to U \to X\oplus U \to Y \to 0$ with $U$ of finite length.
(such a sequence will be called
a {\it Riedtmann-Zwara sequence}). The map $U \to U$ is called
a corresponding {\it steering} map. (Note that in case we deal with modules over
a finite dimensional $k$-algebra and $k$ is an algebraically closed field, then
this notion of degeneration coincides with the usual one, as Zwara [Z2] has shown.)
	\bigskip
The proof of the following result is essentially due to Zwara, he used this
argument in order to show that $Y$ is a degeneration of $X$ if and only if there
is an exact sequence $0 \to Y \to X\oplus V \to V \to 0$ (a co-Riedtmann-Zwara sequence)
with $V$ of finite length.
	\medskip
{\bf Proposition.} {\it Let $X, Y$ be $\Lambda$-modules of finite length.
The following conditions are equivalent:
\item{\rm (1)} $Y$ is a degeneration of $X$.
\item{\rm (2)} There is a Pr\"ufer module $Y[\infty]$ and some natural number $t_0$ such
 that $Y[t+1] \simeq Y[t]\oplus X$ for all $t \ge t_0.$
\item{\rm (3)} There is a Pr\"ufer module $Y[\infty]$ and some natural number $t_0$ such
 that $Y[t_0+1] \simeq Y[t_0]\oplus X$.\par}
	\bigskip
Here is the recipe how to obtain a Pr\"ufer module $Y[\infty]$ starting from a degeneration:
If $Y$ is a degeneration of $X$, say with steering module $U$,
then there exists a monomorphism
$\mu\:U \to U\oplus X$ with cokernel $Y$. The Pr\"ufer module
$Y[\infty]$ we are looking for is 
$$
 Y[\infty] = P(\mu,\left[\smallmatrix 1 \cr 0 \endsmallmatrix\right]).
$$
	\bigskip
Proof of the implication (3) $\implies$ (1).
Assume that there is 
a Pr\"ufer module $Y[\infty]$ such that $Y[t+1] \simeq Y[t]\oplus X.$
We get the following two exact sequences
$$
\gather 
 0 \to Y[t]  \to Y[t+1] \to Y[1] \to 0, \cr 
 0 \to Y[1]  \to Y[t+1] \to Y[t] \to 0, 
\endgather
$$
in the first, the map $Y[t+1] \to Y[1]$ is given by applying $\psi^t$, in the second
the map $Y[t+1] \to Y[t]$ is given by applying $\psi.$ 
In both sequences, we can replace $Y[t+1]$ by $Y[t]\oplus X.$ Thus we obtain 
as first sequence a new Riedtmann-Zwara sequence, and as second sequence a
dual Riedtmann-Zwara sequence:
$$
\gather
 0 \to Y[t] \to Y[t]\oplus X \to Y \to 0, \cr 
 0 \to Y \to Y[t]\oplus X \to Y[t] \to 0,
\endgather
$$
note that both use the same steering module, namely $Y[t].$ Thus:
	\medskip
{\bf Remark.} We see: 
The module $Y$ is a degeneration of $X$ if and only if there 
exists a module $V$ and an exact sequence $0 \to Y \to V \oplus X \to V \to 0$. 
	\medskip
%============================================================
Proof of the proposition. We need further properties of exact squares:
\item{(E2)} {\it For any map $a\:U \to V$, and any module $X$, the following
diagram is exact:}
$$
\CD
 U @>a>> V \cr
 @V{\left[\smallmatrix 1 \cr 0 \endsmallmatrix\right]}VV  
   @VV{\left[\smallmatrix 1 \cr 0 \endsmallmatrix\right]}V \cr
 U\oplus X @>>{a \oplus 1_X}> V\oplus X.
\endCD
$$

\item{(E3)} {\it Let $$
\CD X @>f>> Y_1 \cr
  @V0VV @VVV \cr
  Y_2 @>>f'> Z
\endCD
$$
be exact. Then $f'$ is split mono.}

\item{(E4)} {\it Assume that we have the following exact square
$$
\CD
 U @>a>> V \cr
 @Vb VV  
   @Vb' VV \cr
  W @>>a'> X
\endCD
$$
and that $b$ is a split monomorphism, then the sequence
$$
 0 \to U @>{\left[\smallmatrix a \cr b \endsmallmatrix\right]}>> 
 V \oplus W @>{[b'\ a']}>> X \to 0
$$
splits.}
 	\medskip

Proofs. (E2) is obvious. 	
(E3): Since $\left[\smallmatrix f \cr 0\endsmallmatrix\right]$ is injective, $f\:X \to Y_1$ is
injective. Let $Q$ be the cokernel of $f$. 
We obtain the map $f'$ by forming the induced exact sequence of
$0 \to X @>f>> Y_1 \to Q \to 0$, using the zero map $X \to Y_1$. But such an induced
exact sequence splits.
(E4) Assume that $pb = 1_U.$ Then $[0\ p]\left[\smallmatrix a \cr b \endsmallmatrix\right]
 = 1_U.$
	\bigskip
There is the following lemma  (again, see Zwara [Z1]):
	\medskip
{\bf Lemma (Existence of nilpotent steering maps.)} 
{\it If there is an exact sequence $0 \to U \to X\oplus U \to Y
\to 0,$ then there is an exact sequence $0 \to U' \to X\oplus U' \to Y
\to 0$ such that the map $U' \to U'$ is nilpotent.}
	\medskip
Proof: We can decompose $U = U_1 \oplus U_2 = U_1' \oplus U_2'$ such that
the given map $f\:U \to U$ maps $U_1$ into $U_1'$, $U_2$ into $U_2'$ and such that
the induced maps $f_1\:U_1 \to U_1'$ belongs to the radical of the category,
whereas the induced map $f_2 \:U_2 \to U_2'$ is an isomorphism. We obtain the 
following pair of exact squares 
$$
\CD
 U_1 @>{\left[\smallmatrix 1 \cr 0\endsmallmatrix\right]}>> U_1\oplus U_2 @>>> X \cr
 @Vf_1VV                  @Vf_1\oplus f_2 VV      @VVV  \cr
 U'_1 @>>{\left[\smallmatrix 1 \cr 0\endsmallmatrix\right]}> U'_1\oplus U'_2 @>>> Y 
\endCD
$$
(the left square is exact according to (E2)).
The composition of the squares 
is the desired exact square (note that $U_1'$ is isomorphic to $U_1$).
	\medskip
Assume that a monomorphism $w = \bmatrix \phi \cr g\endbmatrix\:U \to U\oplus X$ 
with cokernel $Y$ and $\phi^t = 0$ is given.
Consider also the canonical embedding
$v= \bmatrix 1 \cr 0\endbmatrix\:U \to U\oplus X$ and form the ladder $U_i(w,v)$
for this pair
of monomorphisms $w,v$. The modules
$Y[i] = U_i(w,v)/U_0(w,v)$ are just the modules we are looking for:
As we know, there is a Pr\"ufer module $(Y[\infty],\psi)$ with $Y[i]$ being the
kernel of $\psi^i.$
	\medskip
We construct the maps $w_n, v_n$ explicitly as follows:
$$
 w_n = \bmatrix \phi \cr
                g \cr
                  & 1_{X^n} \endbmatrix 
 = \bmatrix w \cr
                  & 1_{X^n} \endbmatrix \: U\oplus X^n \to (U\oplus X) \oplus X^n
$$
and 
$$
 v_n = \bmatrix 1_{U\oplus X_n} \cr
                0 \endbmatrix \: U\oplus X^n \to U\oplus X^n \oplus X,
$$
using the recipe (E2). Thus we obtain the following sequence of exact squares:
$$
\CD
 U @>{\left[\smallmatrix \phi \cr g \endsmallmatrix\right]}>> 
 U\!\oplus\! X  @>{\left[\smallmatrix \phi \cr g \cr & 1 \endsmallmatrix\right]}>> 
 U\!\oplus\! X\!\oplus\! X 
 @>{\left[\smallmatrix \phi \cr g \cr & 1 \cr & & 1
 \endsmallmatrix\right]}>> 
 U\!\oplus\! X\!\oplus\! X \!\oplus\! X  &\ \ \longrightarrow \cr
 @V{\left[\smallmatrix 1 \cr 0 \endsmallmatrix\right]}VV
 @V{\left[\smallmatrix 1 \cr & 1 \cr 0 & 0 \endsmallmatrix\right]}VV
 @V{\left[\smallmatrix 1 \cr & 1 \cr & & 1 \cr 0 & 0 & 0 \endsmallmatrix\right]}VV
@V{\left[\smallmatrix 1 \cr & 1 \cr & & 1 \cr & & & 1 \cr
0 & 0 & 0 & 0 \endsmallmatrix\right]}VV \cr
 U\!\oplus\! X  @>>{\left[\smallmatrix \phi \cr g \cr & 1 \endsmallmatrix\right]}>
 U\!\oplus\! X\!\oplus\! X 
 @>>{\left[\smallmatrix \phi \cr g \cr & 1 \cr & & 1
 \endsmallmatrix\right]}>
 U\!\oplus\! X\!\oplus\! X \!\oplus\! X 
 @>>{\left[\smallmatrix \phi \cr g \cr & 1 \cr & & 1 \cr & & & 1 
 \endsmallmatrix\right]}>
 U\!\oplus\! X\!\oplus\! X \!\oplus\! X \!\oplus\! X &\ \ \longrightarrow \cr
\endCD
$$
In particular, we have $U_n = U_n(w,v) = U\oplus X^n.$

Note that the composition $w_{n-1}\cdots w_0 \:U \to U\oplus X^n$ is of the form
$\bmatrix \phi^n \cr g_n \endbmatrix$ for some $g_n\:U \to X^n$.

We also have the following sequence of exact squares:
$$
\CD 
 U=U_0  @>w_0>> 
 U_1  @>w_1>> 
 U_2  @>w_2>> 
 U_3  @>w_3>> \cdots \cr
@VVV 
@VVV 
@VVV 
@VVV  \cr
 0  @>>> 
 Y[1]  @>s_1>> 
 Y[2]  @>s_2>> 
 Y[3]  @>s_3>> \cdots
\endCD
$$
where the vertical maps are of the form
$$
  U_n = U\oplus X^n @>{[ h_n \ q_n ]}>> Y[n].
$$

The composition of these exact squares yields an exact square
$$
\CD U @>{w_{n-1}\cdots w_0}>> U\oplus X^n \cr
 @VVV    @VV{[ h_n \ q_n ]}V \cr
 0 @>>> Y[n]
\endCD
$$
Here we may insert the following observation: This sequence shows that {\it
the module $Y[n]$ is a degeneration of the module $X^n$.}

Since the composition $w_{n-1}\cdots w_0 \:U \to U\oplus X^n$ is of the form
$\bmatrix \phi^n \cr g_n \endbmatrix,$ and $\phi^t = 0,$ 
it follows that $h_t$ is a split 
monomorphism, see (E3). 

Also, we can consider the following two exact squares, with 
$w = \bmatrix \phi \cr g\endbmatrix\:U \to V = U\oplus X$ (the upper square
is exact, according to (E2)):
$$
\CD 
 U @>w>> V \cr
 @V{\left[\smallmatrix 1 \cr 0 \endsmallmatrix\right]}VV 
 @VV{\left[\smallmatrix 1 \cr 0 \endsmallmatrix\right]}V \cr
 U\oplus X^t @>{\left[\smallmatrix w \cr & 1 \endsmallmatrix\right]}>> V\oplus X^t \cr
 @V{[ h_t \ q_t ]}VV   @VV{[ h_{t+1} \ q_{t+1} ]}V \cr
 Y[t] @>>> Y[t\!+\!1]
\endCD
$$
The vertical composition on the left is $h_t$, thus, as we have shown, a split monomorphism.
This shows that the exact sequence corresponding to the composed square splits (E4):
This yields
$$
 U \oplus Y[t\!+\!1] \simeq Y[t] \oplus V = Y[t] \oplus U \oplus X.
$$
Cancelation of $U$ gives the desired isomorphism:
$$
 Y[t\!+\!1] \simeq  Y[t] \oplus X.
$$
	\bigskip
{\bf Remark to the proof.} Given the Riedtmann-Zwara sequence
$$
 0 \to U @>{\left[\smallmatrix \phi \cr g \endsmallmatrix\right]}>>
 U \oplus X @>>> Y \to 0,
$$
we have considered the following pair of monomorphisms
$$
 w = \left[\smallmatrix 1 \cr 0 \endsmallmatrix\right], 
 w' =  \left[\smallmatrix \phi \cr g \endsmallmatrix\right]\: U \to U\oplus X.
$$
The corresponding Pr\"ufer modules are $X^{(\infty)}$ and $Y[\infty]$, respectively. 
And $U_n(w,w') = U\oplus X^n.$ As we know, we can assume that $\phi$ is nilpotent.
Then all the linear combinations 
$$
 w+\lambda w' = \left[\smallmatrix 1+\lambda \phi \cr g \endsmallmatrix\right]
$$
with $\lambda \in k$
are also split monomorphisms (with retraction 
$\left[\smallmatrix \eta & 0 \endsmallmatrix\right]$, where 
$\eta = (1+\lambda\phi)^{-1}$).
	\bigskip
{\bf Corollary.} {\it Assume that $Y$ is a degeneration of $X$. Then there exists
a Pr\"ufer module $Y[\infty]$ such that $Y[\infty]$ is isomorphic to $Y[t]\oplus X^{(\omega)}$
for some natural number $t$.}

	\bigskip\bigskip
%=====================================================================
{\bf 5. Application: The theorem of Bautista-Perez.}
	\medskip
Here we assume that we deal with an artin algebra $\Lambda$, and all the modules
are $\Lambda$-modules of  finite length. 
	\medskip
{\bf Proposition.} {\it Let $W$ be a module with $\Ext^1(W,W) = 0$ and assume there is
given an exact sequence $0 \to U \to V \to W \to 0$. Then the cokernel of any monomorphism
$U \to V$ is a degeneration of $W$.}
	
	\bigskip
{\bf Corollary (Bautista-Perez).} {\it Let $U, V$ be modules, and let 
$W$ and $W'$ be cokernels of monomorphisms $U \to V.$
Assume that both $\Ext^1(W,W) = 0$ and
$\Ext^1(W',W') = 0.$ Then the modules $W$ and $W'$ are isomorphic.}
	\medskip 
Both assertions are well-known in case $k$ is an algebraically closed field: in this
case, the conclusion of proposition just asserts that $W'$ is a degeneration 
of $W$ in the sense of
algebraic geometry.
The main point here is to deal with
the general case when $\Lambda$ is an arbitrary artin algebra. 
The corollary stated above (under the additional assumptions that $V$ is projective
and that $w(U), w'(U)$ are contained in the radical of $V$) is due
to Bautista and Perez [BP] and this result was presented 
by Smal{\o} with a new proof [S]
at Mar del Plata. 
	\medskip
We need the following well-known lemma.
	\medskip
{\bf Lemma.} {\it Let $W$ be a module with $\Ext^1(W,W) = 0.$
Let $U_0 \subset U_1 \subset U_2 \subset \cdots$ be a sequence of inclusions
of modules with $U_{i}/U_{i-1}  = W$ for all $i \ge 1.$ Then there is a natural
number $n_0$ such that $U_{n} \subset U_{n+1}$ is a split monomorphism for all $n \ge n_0.$}
	\medskip
Let us use it in order to finish the proof of proposition.
Let $U_0 = U, U_1 = V,$ and $w_0\:U_0 \to V_0$ the given monomorphism with cokernel $W$.
Let $v_0\:U_0\to U_1$ be an additional monomorphism, say with cokernel $W'$. 
Thus we are in the setting of section 1.
We apply Lemma to the chain of inclusions
$$
U_0 @>w_0>> U_1 @>w_1>> U_2 @>w_2>> \cdots
$$
and see that there is $n$ such that $w_{n}\:U_{n} \to U_{n+1}$ splits. This shows that
$U_{n+1}$ is isomorphic to $U_{n}\oplus W.$ But we also have the exact sequence
$$
 0 \to U_{n} @>v_n>> U_{n+1} \to W' \to 0.
$$
Replacing $U_{n+1}$ by $U_{n}\oplus W$, we see that we get an exact sequence of the form
$$
 0 \to U_{n} @>v_n>> U_{n}\oplus W \to W' \to 0
$$
(a Riedtmann-Zwara sequence), as asserted.
	\medskip
Proof of Corollary. It is well-known that the existence of exact sequences
$$
 0 \to X \to X\oplus W \to W' \to 0  \t{and}  0 \to Y \to Y\oplus W' \to W \to 0
$$
implies that the modules $W$ and $W'$ are isomorphic [Z1].
But in our case we just have to 
change  one line in the proof of proposition  in order to get the required
isomorphism. Thus, assume that both $\Ext^1(W,W) = 0$ and 
$\Ext^1(W',W') = 0$. Choose $n$ such that {\bf both} the inclusion maps
$$
  w_{n}\:U_{n} \to U_{n+1} \t{and} v_{n}\:U_{n} \to U_{n+1}
$$
split. Then $U_{n+1}$ is isomorphic both to 
$U_{n}\oplus W$ and to $U_{n}\oplus W'$, thus it follows from the Krull-Remak-Schmidt
theorem that $W$ and $W'$ are isomorphic.
	\medskip
{\bf Remark.} {\it Assume that $w, w'\:U, V$ are monomorphisms with cokernels 
$W$ and $W'$, respectively, and that $\Ext^1(W,W) = 0$ and
$\Ext^1(W',W') = 0.$ Then $w$ splits if and only if $w'$ splits.}
	\medskip
Proof: According to the corollary, we can assume $W = W'$. Assume that $w$ splits, thus
$V$ is isomorphic to $U\oplus W$. Look at the exact sequence $0 \to U @>w'>> V \to W \to 0$.
If it does not split, then $\dim\End(V) < \dim\End(U\oplus W)$, but $V$ is isomorphic to $U\oplus 
W.$
	\bigskip
As we have mentioned, the lemma is well-known; an equivalent assertion 
was used for example
by Roiter in his proof of the first Brauer-Thrall conjecture, a corresponding proof can be
found in [R4]. We include here a slightly different proof:

Applying the functor $\Hom(W,-)$ to the short exact sequence 
$0 \to U_{i-1} @>w_{i-1}>> U_i \to W \to 0,$
we obtain the exact sequence
$$
 \Ext^1(W,U_{i-1}) \to \Ext^1(W,U_i) \to \Ext^1(W,W).
$$
Since the latter term is zero, we see that we have a sequence of surjective maps
$$
 \Ext^1(W,U_0) \to \Ext^1(W,U_1) \to \cdots \to \Ext^1(W,U_i) \to \cdots,
$$
being induced by the inclusion maps $U_0 \to U_1 \to \cdots \to U_i \to \cdots.$
The maps between the $\Ext$-groups 
are $k$-linear. Since $\Ext^1(W,U_0)$ is a $k$-module of finite length, the
sequence of surjective maps must stabilize: there is some $n_0$ such that the inclusion
$U_{n} \to U_{n+1}$ induces an isomorphism 
$$
 \Ext^1(W,U_{n}) \to \Ext^1(W,U_{n+1}) 
$$
for all $n \ge n_0.$ Now we consider also some $\Hom$-terms: the exactness of
$$
 \Hom(W,U_{n+1}) \to \Hom(W,W) \to \Ext^1(W,U_{n}) \to \Ext^1(W,U_{n+1})
$$
shows that the connecting homomorphism is zero, and thus that the map $\Hom(W,U_{n+1}) \to \Hom(W,W)$
(induced by the projection map $p\:U_{n+1} \to W$) is surjective. But this means that there
is a map $h\in \Hom(W,U_{n+1})$ with $ph = 1_W$, thus $p\: U_{n+1} \to W$ is a split
epimorphism and therefore the inclusion map $U_{n} \to U_{n+1}$ is a split monomorphism. 
	\medskip
{\bf Remark.} In general, there is no actual bound on the number $n_0$.
However, in case of dealing with the chain of inclusions
$$
  U_0 @>w_0>> U_{1} @>w_1>> U_2 @>w_n>> \cdots
$$
such a bound exists, namely the length of
$\Ext^1(W,U_0)$ as a $k$-module, or, even better, the length of
$\Ext^1(W,U_0)$ as an $E$-module, where $E = \End(W).$
	\medskip
Proof: Look at the surjective maps
$$
 \Ext^1(W,U_0) \to \Ext^1(W,U_1) \to \cdots \to \Ext^1(W,U_i) \to \cdots,
$$
being induced by the  maps $U_n @>w_n>> U_{n+1}$ (and these maps are not only
$k$-linear, but even $E$-linear). 
Assume that $ \Ext^1(W,U_{n}) \to \Ext^1(W,U_{n+1})$ is bijective, for some $n$.
As we have seen above, this implies that the sequence
$$
 0 \to U_n @>w_n>> U_{n+1} \to W \to 0 \tag{$*$}
$$ 
splits. Now the map $w_{n+1}$ is obtained from $(*)$ as the induced exact sequence 
using the map $w'_n$. With $(*)$ also any induced exact sequence will split.
Thus $w_{n+1}$ is a split monomorphism (and  
$ \Ext^1(W,U_{n+1}) \to \Ext^1(W,U_{n+2})$ will be bijective, again). Thus, as soon
as we get a bijection $ \Ext^1(W,U_{n}) \to \Ext^1(W,U_{n+1})$ for some $n$, then also
all the following maps $ \Ext^1(W,U_{m}) \to \Ext^1(W,U_{m+1})$ with $m > n$ are 
bijective.

	\bigskip
{\bf Example.} Consider the $D_4$-quiver with subspace orientation: 
$$
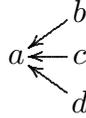
{\beginpicture
\setcoordinatesystem units <.7cm,.5cm>
\put{$a$} at 0 1
\put{$b$} at 1 2
\put{$c$} at 1 1
\put{$d$} at 1 0
\arr{0.8 0.2}{0.2 0.8}
\arr{0.8 1}{0.2 1}
\arr{0.8 01.8}{0.2 1.2}
\endpicture}
$$
and let $\Lambda$ be its path algebra over some field $k$. 
We denote the indecomposable $\Lambda$-modules by the corresponding dimension vectors.
Let 
$$U_0 = \smallmatrix   &0 \cr
                      1&0 \cr
                       &0 \endsmallmatrix,\quad
U_1 = \smallmatrix     &1 \cr
                      2&1 \cr
                       &1 \endsmallmatrix,\quad 
W = \smallmatrix       &1 \cr
                      1&1 \cr
                       &1 \endsmallmatrix,\quad 
W' = \smallmatrix      &0 \cr
                      1&1 \cr
                       &1 \endsmallmatrix \oplus 
         \smallmatrix  &1 \cr
                      0&0 \cr
                       &0 \endsmallmatrix.
$$
Note that a map $w_0\:U_0 \to U_1$ with cokernel $W$ exists only in case the base-field $k$
has at least 3 elements; of course, there is always a map $w'_0\:U_0 \to U_1$ with cokernel $W'.$

We have $\dim\Ext^1(W,U_0) = 2,$ and it turns out that the module $U_2$ is the following:
$$
 U_2 =   \smallmatrix  &0 \cr
                      1&1 \cr
                       &1 \endsmallmatrix  \oplus
         \smallmatrix  &1 \cr
                      1&0 \cr
                       &1 \endsmallmatrix  \oplus
         \smallmatrix  &1 \cr
                      1&1 \cr
                       &0 \endsmallmatrix.
$$

The pushout diagram involving the modules $U_0,\ U_1$ (twice) and $U_2$ is
constructed as follows: denote by $\mu_a, \mu_b, \mu_c$ monomorphisms $U_0 \to U_1$
which factor through the indecomposable projective modules $P(a), P(b), P(c)$,
respectively. We can assume that $\mu_c = -\mu_a-\mu_b$, so that a mesh relation
is satisfied. Denote the 3 summands of $U_2$ by $M_a, M_b, M_c$, with non-zero maps
$\nu_a\:U_1 \to M_a,\ \nu_b\:U_1 \to M_b,\ \nu_c\:U_1 \to M_c,$ such that
$\nu_a\mu_a = 0,\ \nu_b\mu_b = 0,\ \nu_c\mu_c = 0.$   There is the following
commutative square, for any $q\in k,$ we are interested when $q\notin\{0,1\}$:
$$
\CD
 U_0 @>{w_0 = \mu_a+q\mu_b}>> U_1 \cr
 @V{v_0=\mu_a}VV           @VV{v_1=\left[\smallmatrix 0 \cr \nu_b \cr \nu_c\endsmallmatrix\right]}V \cr
 U_1 @>>{w_1 = \left[\smallmatrix \nu_a \cr \nu_b \cr (1-q)\nu_c\endsmallmatrix\right]}> U_2 \cr
\endCD
$$
(the only calculation which has to be done concerns the third entries: 
$\nu_c(\mu_a+q\mu_b) = (1-q)\nu_c\mu_a$). 
Note that $w_1$ (as well as $w'_1$) does not split.

But now we deal with a module $U_2$ such that $\Ext^1(W,U_2) = 0.$ This implies that $U_3$
is isomorphic to $U_2\oplus W$. Thus the next pushout construction yields an exact sequence
of the form
$$
 0 \to U_2 \to U_2\oplus W \to W' \to 0.
$$
	\bigskip\bigskip
%==========================================================================
{\bf References.}
	\medskip
\item{[BP]} Bautista, R. and Perez, E.: On modules and complexes without self-extensions.
   Communications in Algebra 34, 3139-3152.
\item{[R1]} Ringel, C.M.: Representations of k-species and bimodules. J.Algebra 
 41 (1976), 269-302. 
\item{[R2]} Ringel, C.M.: Bricks in hereditary length categories.
  Resultate der Math. 6 (1983), 64-70. 

\item{[R3]} Ringel, C.M.: Some algebraically compact modules I. 
In: Abelian Groups and Modules (ed. A. Facchini and C. Menini). Kluwer (1995), 419-439. 

\item{[R4]} Ringel, C.M.: The Gabriel-Roiter measure. Bull. Sci. math. 129 (2005). 726-748.
\item{[R5]} Ringel, C.M.: Pr\"ufer modules which are not of finite type. In preparation.
\item{[S]} Smal{\o}, S.: Lectures on Algebras. Mar del Plata, Argentina March  2006.
   To appear in Revista Uni\'on Matem\'atica Argentina. 
\item{[Z1]} Zwara, G.: A degeneration-like order for modules. Arch\. Math\. 71 (1998), 437-444
\item{[Z2]} Zwara, G.: Degenerations of finite-dimensional modules are given
   by extensions. Compositio Mathematica 121 (2001), 205-218.

\bye